\documentclass[a4paper,10pt,oneside]{article}   
\usepackage[english]{babel}
\usepackage[utf8]{inputenc}
\usepackage[T1]{fontenc}
\usepackage{amssymb,amsthm,amsmath}
\usepackage{verbatim}
\usepackage{enumitem} 
\usepackage{a4wide}
\usepackage{hyperref}
\usepackage{definitions}

\usepackage{graphicx}

\usepackage[sorting=nyt, style=alphabetic, firstinits=true, doi=false,isbn=false,url=false, backend=biber]{biblatex}
\bibliography{basic}
\AtEveryBibitem{\ifentrytype{book}{\clearfield{pages}}}

\usepackage{xcolor}
\hypersetup{
    colorlinks,
    linkcolor={red!50!black},
    citecolor={blue!50!black},
    urlcolor={blue!80!black}
}

\begin{document}

\title{The Goldbach conjecture with summands in arithmetic progressions}
\author{Juho Salmensuu}
%\address{Department of Mathematics and Statistics\\ University of Turku\\ FI-20014 University of Turku\\ Finland}
%\email{juelsa@utu.fi}
\date{}

\maketitle

\begin{abstract}
   We prove that, for almost all $r \leq N^{1/2}/\log^{O(1)}N$, for any given $b_1 \Mod r$ with $(b_1, r) = 1$, and for almost all $b_2 \Mod r$ with $(b_2, r) = 1$, we have that almost all natural numbers $2n \leq N$ with $2n \equiv b_1 + b_2 \Mod r$ can be written as the sum of two prime numbers $2n = p_1 + p_2$, where $p_1 \equiv b_1 \Mod r$ and $p_2 \equiv b_2 \Mod r$. This improves the previous result which required $r \leq N^{1/3}/\log^{O(1)}N$ instead of $r \leq N^{1/2}/\log^{O(1)}N$. We also improve some other results concerning variations of the problem.
\end{abstract}

\section{Introduction}

\subsection{Statements of results}
In this paper, we study the binary Goldbach conjecture in the case where the summands are restricted to arithmetic progressions with large moduli. Let $N > 0$ and define
\begin{displaymath} 
E_{b_1, b_2, r}(N) := \# \{2n \leq N \mid 2n \equiv b_1 + b_2 \Mod r,  2n \not= p_1 + p_2 \text{ for all primes } p_i \equiv b_i \Mod r\}.
\end{displaymath}

Our main result is the following.

\begin{Theorem} \label{theorem_main} Let $A, B> 0$, $N \geq 3$ and $R \leq N^{1/2} / \log^{A} N$. Then, for all but $O(R / \log^{B} N)$ integers $3 \leq r \leq R$, for any fixed $b_1 \Mod r$ with $(b_1, r) = 1$ and for all but $O(r / \log^{B} r)$ integers $b_2 \Mod r$ with $(b_2, r) = 1$, we have 
\begin{displaymath}
E_{b_1, b_2, r}(N) = o(N/r),
\end{displaymath} 
provided that $A$ is large enough depending on $B$. 
\end{Theorem}

This improves Bauer's result \cite{bauer_latest}, which had $1/3$ instead of $1/2$ and had restriction that $r$ has to be a prime. When also $b_2$ is fixed, Bauer \cite{bauer_large_sieve} has showed that $E_{b_1, b_2, r}(N) = o(N/r)$ holds for almost all primes $r \leq N^{7/30 - \epsilon} / \log^{O(1)} N$. 

We will also prove the following theorem, where the exceptional set is very small.

\begin{Theorem} \label{theorem_3} Let $\epsilon> 0$ and $N \geq 3$. There exists $D > 0$ such that, for all but $O(\log^{D} N)$ primes $r \leq N^{5/24 - \epsilon}$, for any fixed $b_1, b_2 \Mod r$ with $(b_1, r) = (b_2, r) = 1$, we have 
\begin{displaymath}
E_{b_1, b_2, r}(N) = o(N/r). 
\end{displaymath} 
\end{Theorem}

This improves Bauer's result \cite{bauer_small_exceptional_set}, which had $5/48$ instead of $5/24$. \cite[Theorem 2]{bauer_latest} has $5/24$, but the minor arc calculations appear to be incorrect. Using the minor arcs calculations from \cite{bauer_old} one can recover \cite[Theorem 2]{bauer_latest} in the following form:
\\\\
\noindent \textbf{Theorem A.} \textit{Let $B, \epsilon> 0$ and $N \geq 3$. There exists $D > 0$ such that, for all but $O(\log^{D} N)$ primes $3 \leq r \leq N^{5/24 - \epsilon}$, for any fixed $b_1 \Mod r$ with $(b_1, r) = 1$ and for all but $O(r / \log^{B} r)$ integers $b_2 \Mod r$ with $(b_2, r) = 1$, we have}
\begin{displaymath}
E_{b_1, b_2, r}(N) = o(N/r). 
\end{displaymath}

We improve this result in the following theorem. 

\begin{Theorem} \label{theorem_2} Let $ B, \epsilon> 0$ and $N \geq 3$. There exists $D > 0$ such that, for all but $O(\log^{D} N)$ primes $3 \leq r \leq N^{5/12 - \epsilon}$, for any fixed $b_1 \Mod r$ with $(b_1, r) = 1$ and for all but $O(r / \log^{B} r)$ integers $b_2 \Mod r$ with $(b_2, r) = 1$, we have 
\begin{displaymath}
E_{b_1, b_2, r}(N) = o(N/r). 
\end{displaymath} 
\end{Theorem}

Assuming GRH, we can completely dispose of the exceptional sets for the moduli from Theorems \ref{theorem_main}, \ref{theorem_3} and \ref{theorem_2} and replace the exponent $5/24$ in Theorem \ref{theorem_3} with $1/4$.  

\begin{center}
\begin{tabular}{ c || c | c || c | c || c | c || c }
 & 
 \rotatebox{90}{Old} &  
 \rotatebox{90}{New} & 
 \rotatebox{90}{Old} &
 \rotatebox{90}{New}&
 \rotatebox{90}{Old} & 
 \rotatebox{90}{New} & 
 \rotatebox{90}{Latest} \\ \hline
Restricted to primes  & x   & -    & x     &  x    & x    & x    & x    \\ 
Large exceptional set & x   & x    & -     &  -    & -    & -    & x    \\
Small exceptional set & -   & -    & x     &  x    & x    & x    & -    \\
Mean value over $b_1$ & x   & x    & -     &  -    & x    & x    & -     \\ \hline
Exponent 			  & 1/3 & 1/2  & 5/48  &  5/24 & 5/24 & 5/12 & 7/30   \\        
\end{tabular}
\\
The latest results and the improvements.
\end{center}

\textbf{Remark.} Using the methods of this paper it is possible to prove similar results for the ternary Goldbach problem. For example, we can prove the following: Let $A, B> 0$, $N > 3$ and $R \leq N^{1/2} / \log^{A} N$. Then, for all but $O(R / \log^{B} N)$ integers $3 < r \leq R$, for any fixed $b_1, b_2 \Mod r$ with $(b_1b_2, r) = 1$, for all but $O(r / \log^{B} r)$ integers $b_3 \Mod r$ with $(b_3, r) = 1$ and for all $n \leq N$ with $n \equiv b_1 + b_2 + b_3 \Mod r$, we can write
\begin{displaymath}
n = p_1 + p_2 + p_3,
\end{displaymath} 
where $p_1, p_2, p_3$ are primes with $p_i \equiv b_i \Mod r$.

\subsection{Notation}

Denote divisor function by $\tau(n) := \sums{d | n} 1$, Euler's totient function by 
\begin{displaymath}
\phi(n) := \sums{t \leq n \\ (t, n) = 1} 1
\end{displaymath}
and von Mangoldt function by
\begin{displaymath}
\Lambda(n) :=
\begin{cases}
\log p & \text{ if } n = p^k, \\
0 & \text{ otherwise.}
\end{cases}
\end{displaymath}
Set
\begin{equation} \label{definition_of_E}
E_{d}(x) := \max_{t \leq x} \max_{h: (h, d) = 1} \Big| \sums{n \leq t\\ n \equiv h \Mod{d}} \Lambda(n) - \frac{t}{\phi(d)} \Big| + 1.
\end{equation}

We use abbreviations $e(x) := e^{2\pi i x}$ and $e_q(n) := e(n/q)$.

\indent Let $f: \R \rightarrow \C$ and $g: \R \rightarrow \R_+$. We write
$f = O(g), f \ll g$
if there exists a constant $C > 0$ such that $|f(x)| \leq C g(x)$
for all values of $x$ in the domain of $f$. If the implied constant $C$ depends on some contant $\epsilon$ we use notations $O_\epsilon, \ll_\epsilon$. We also write $f = o(g)$ if
\begin{displaymath}
\lim _{x \rightarrow \infty} \frac{f(x)}{g(x)} = 0.
\end{displaymath}

Notation $n \sim N$ means $N \leq n < 2N$. For $\alpha \in \R$ denote $||\alpha|| = \min_{n \in \Z} |\alpha - n|$.

\subsection{Outline}

In this section, we present our main ideas used to prove the theorems. 

Let 
\begin{displaymath}
\mathfrak{S}_r(h) := \prod_{\substack{p \nmid r\\ p \nmid h}} \Big( 1 - \frac{1}{(p-1)^2}\Big) \prod_{\substack{p \nmid r\\ p | h}} \Big( 1 + \frac{1}{(p-1)}\Big)
\end{displaymath}
be the singular series for Goldbach's problem in arithmetic progressions.

Our aim is to prove the following three theorems. 

\begin{Theorem} \label{theorem_technical_main} Let $A > 0$. Let $R> 0$ and $N \geq 3$ be such that $R \leq (RN)^{1/2} / \log^{2A} N$. Then
\begin{displaymath}
\sums{r \leq R} \max_{b_1: (b_1, r) = 1} \sums{b_2 \Mod r\\ (b_2, r) = 1} \sums{n \leq N} \Big| \sums{n_1, n_2\\ n = n_1 + n_2} \Lambda(rn_1 + b_1)\Lambda(rn_2 + b_2) - \frac{r^2}{\phi(r)^2}\mathfrak{S}_r(rn + b_1 +b_2)n \Big | \ll_{B} \frac{N^2 R^2}{\log^B N},
\end{displaymath}
for any $B > 0$, provided that $A$ is large enough depending on $B$.
\end{Theorem}

\begin{Theorem} \label{theorem_technical_3} Let $\epsilon, B > 0$ and $N \geq 3$. There exists $D = D(B) > 0$ such that, for all but $O(\log^D N)$ primes $r$ with $r \leq (rN)^{5/24 - \epsilon}$, we have
\begin{displaymath}
\max_{b_1, b_2: (b_1b_2, r) = 1} \sums{n \leq N} \Big| \sums{n_1, n_2\\n = n_1 + n_2} \Lambda(rn_1 + b_1)\Lambda(rn_2 + b_2) - \frac{r^2}{\phi(r)^2}\mathfrak{S}_r(rn + b_1 +b_2)n \Big | \ll_{B, \epsilon} \frac{N^2}{\log^B N}.
\end{displaymath}
\end{Theorem}

\begin{Theorem} \label{theorem_technical_2} Let $\epsilon, B > 0$ and $N \geq 3$. There exists $D = D(B) > 0$ such that, for all but $O(\log^D N)$ primes $r$ with $r \leq (rN)^{5/12 - \epsilon}$, we have
\begin{displaymath}
\max_{b_1: (b_1, r) = 1} \sums{b_2 \Mod r\\ (b_2, r) = 1} \sums{n \leq N} \Big| \sums{n_1, n_2\\n = n_1 + n_2} \Lambda(rn_1 + b_1)\Lambda(rn_2 + b_2) - \frac{r^2}{\phi(r)^2}\mathfrak{S}_r(rn + b_1 +b_2)n \Big | \ll_{B, \epsilon} \frac{N^2 r}{\log^B N}.
\end{displaymath}
\end{Theorem}

Theorems \ref{theorem_main}, \ref{theorem_3} and \ref{theorem_2} follow respectively from Theorems \ref{theorem_technical_main}, \ref{theorem_technical_3} and   \ref{theorem_technical_2}. (Note that the contribution of the prime powers in the sums is negligible.)

We use the circle method to prove these theorems. The major improvement to the previous results comes from the way we arrange the circle method. Let
\begin{displaymath}
S_{b, r}(H, \alpha) := \sums{n \leq H} \Lambda(rn + b) e(\alpha n), \ S_{b, r}'(H, \alpha) := \sums{n \leq H\\ n \equiv b \Mod r} \Lambda(n) e(\alpha n)
\end{displaymath}
and $M = rN + b_1 + b_2$. In the previous papers, the circle method has been applied in the following way:
\begin{displaymath}
\sums{M = n_1 + n_2 \\ n_i \equiv b_i \Mod{r_i}} \Lambda (n_1)\Lambda (n_2)  = \int_0^1 S_{b_1, r}'(M, \alpha)S_{b_2, r}'(M, \alpha)  e(-\alpha M)  d\alpha.
\end{displaymath} 
We apply the circle method inside the arithmetic progression ($M \equiv b_1 + b_2 \Mod r$):
\begin{equation*} 
\sums{N = n_1 + n_2} \Lambda (rn_1 + b_1)\Lambda (rn_2 + b_2)  = \int_0^1 S_{b_1, r}(N, \alpha)S_{b_2, r}(N, \alpha)  e(-\alpha N)  d\alpha.
\end{equation*} 
Using the circle method inside the arithmetic progression leads to easier exponential sums. We demonstrate this by comparing Type I estimates in both cases. 

Let $\alpha \in [0, 1]$ and $a, q \in \N$ be such that $1 \leq a \leq q$, $(a, q) = 1$ and $|\alpha - a/q| \leq q^{-2}$. Let $N, M \geq 1$.  Let $a_n$ be a complex sequence such that $|a_n| \leq 1$. Write $X := NM$.

Using the standard methods to evaluate type I sum corresponding to $S'$, we see that
\begin{equation} \label{eq_type_I_first}
\sums{n \sim N, m \sim M\\ nm \equiv b \Mod r} a_m e(\alpha nm) \ll \frac{X}{r \log^A X},
\end{equation} 
provided that $r \log^{A'} X \leq q \leq X / (r \log^{A'} X)$ and $M \leq X / (r^2 \log^{A'} X)$ for some $A' > 0$ depending on $A$. Similarly, for type I sum corresponding to $S$, we see that
\begin{equation}\label{eq_type_I_second}
\sums{n \sim N, m \sim M\\ nm \equiv b \Mod r} a_m e_r(\alpha nm) \ll \frac{X}{r \log^A X},
\end{equation} 
provided that $\log^{A'} X \leq q \leq X / (r \log^{A'} X)$ and $M \leq X / (r \log^{A'} X)$ for some $A' > 0$ depending on $A$. We can see that (\ref{eq_type_I_second}) holds for a much wider range than (\ref{eq_type_I_first}).

For $A, N \geq 1$ define 
\begin{displaymath}
\mathfrak{M}_{A, N} := \bigcup_{\substack{q \leq \log^{A} N\\ 1 \leq a \leq q\\ (a, q) = 1}} \mathfrak{M}_{A, N}(q, a),
\end{displaymath}
where 
\begin{displaymath}
\mathfrak{M}_{A, N}(q, a) := \{\alpha \in [0, 1] \mid |\alpha - a/q| \leq N^{-1} \log^{A} N\}.
\end{displaymath}
Also define $\mathfrak{m}_{A} := \mathfrak{m}_{A, N} := [0, 1] \setminus \mathfrak{M}_{A, N}$. We call $\mathfrak{M}_{A, N}$ major arcs and $\mathfrak{m}_{A, N}$ minor arcs. We split 
\begin{align*}
\int_0^1 S_{b_1, r}(N, \alpha)S_{b_2, r}(N, \alpha)  e(-\alpha N)  d\alpha = & \int_{\mathfrak{M}_{A, N}} S_{b_1, r}(N, \alpha)S_{b_2, r}(N, \alpha)  e(-\alpha N)  d\alpha \\
& \ + \int_{\mathfrak{m}_{A, N}} S_{b_1, r}(N, \alpha)S_{b_2, r}(N, \alpha)  e(-\alpha N)  d\alpha.
\end{align*}

In Section \ref{section_major_arcs} we give an asymptotic formula for the major arcs contribution. We analyse the minor arc contribution in Section \ref{section_minor_arcs}. Since the error term, coming from the minor arc and major arc analysis, in some cases depends on $E_d(x)$ defined in (\ref{definition_of_E}), we will give some estimates for $E_d(x)$ in Section \ref{section_error_term}. We prove Theorems \ref{theorem_technical_main}, \ref{theorem_technical_3} and  \ref{theorem_technical_2} in Section \ref{section_proofs}.

\section{Auxiliary lemmas}
In this section, we present some auxiliary lemmas, which we use later. 

\begin{Lemma}  \label{lemma_geometric_series} Let $\alpha \in \R$ and $N_1, N_2 \in \N$ with $N_1 < N_2$. Then
\begin{displaymath}
\sums{N_1 < n  \leq N_2}e(\alpha n) \ll \min(N_2 - N_1, ||\alpha||^{-1}).
\end{displaymath}
\end{Lemma}

\begin{proof} See e.g. \cite[Lemma 4.7]{nathanson}.
\end{proof}

\begin{Lemma} \label{lemma_xy_exp_sum} For any $X, Y \geq 1$, we have
\begin{displaymath}
\sums{1 \leq m \leq X} \min(Y, ||\alpha m||^{-1}) \ll \Big( \frac{XY}{q} + X +q \Big) \log 2qX.
\end{displaymath}
\end{Lemma}

\begin{proof} Trivially
\begin{displaymath}
\sums{1 \leq m \leq X} \min(Y, ||\alpha m||^{-1}) \leq \sums{1 \leq m \leq X} \min\Big(\frac{XY}{m}, ||\alpha m||^{-1}\Big)
\end{displaymath}
and the claim follows from a standard estimate (see e.g. \cite[Lemma 4.10]{nathanson}).
\end{proof}

\begin{Lemma} \textbf{(Vaughan's identity)}  For any $y \geq 1$, $n > y$, we have
\begin{displaymath}
\Lambda(n) = \sums{b | n\\ b \leq y} \mu(b) \log \frac{n}{b} - \sums{bc | n \\ b,c \leq y} \mu(b) \Lambda(c) + \sums{bc | n\\ b, c > y} \mu(b) \Lambda(c).
\end{displaymath}
\end{Lemma}

\begin{proof}
See e.g. \cite[Proposition 13.4]{iwaniec-kowalski}.
\end{proof}

\noindent

\section{Major arcs} \label{section_major_arcs}

In this section, our aim is to prove that 
\begin{displaymath}
\int_{\mathfrak{M}_{A, N}} S_{b_1, r}(n, \alpha) S_{b_2, r}(n, \alpha) e(-\alpha n) d\alpha  \approx \frac{r^2}{\phi(r)^2}\mathfrak{S}(rn + b_1 +b_2)n.
\end{displaymath}
We use the standard circle method machinery to do so.

\subsection{Generating function}

In this subsection, we prove an approximation lemma for the generating function
\begin{displaymath}
S_{b, r}(N, \alpha) = \sums{n \leq N} \Lambda(rn + b) e(\alpha n).
\end{displaymath}

First, we prove the following rational exponential sum estimate. 

\begin{Lemma} \label{lemma_rational_exp_sum}  Let $a, b, q, r \in \N$ be such that $(a, q) = (b, r) = 1$. Then
\begin{displaymath}
\sums{h = 1\\ h \equiv b \Mod r\\ (h, q) = 1}^{rq} e_{rq}(ah) =  \mathbf{1}_{(q, r) = 1} \mu(q)e_{rq}(abq^{\phi(r)}).
\end{displaymath}
\end{Lemma}

\begin{proof} Since $(b, r) = 1$, we have for $d | q$ that the congruence system
\begin{align*}
\begin{cases}
x  \equiv b \Mod r \\
x  \equiv 0 \Mod d 
\end{cases}
\end{align*}
is soluble only if $(d, r) = 1$ in which case $x \equiv bd^{\phi(r)}$ is the unique solution $\Mod{dr}$. Hence
\allowdisplaybreaks
\begin{eqnarray*}
\sums{h = 1\\ h \equiv b \Mod r\\ (h, q) = 1}^{rq} e_{rq}(ah) 
&=& \sums{h = 1\\ h \equiv b \Mod r}^{rq} e_{rq}(ah) \sums{d | (h, q)} \mu(d) \\
&=& \sums{d | q \\ (d, r) = 1} \mu(d) \sums{h = 1\\ h \equiv b \Mod r \\ h \equiv 0 \Mod d}^{rq} e_{rq}(ah) \\
&=& \sums{d | q \\ (d, r) = 1} \mu(d) \sums{h = 1\\ h \equiv bd^{\phi(r)} \Mod{rd} }^{rq} e_{rq}(ah) \\
&=& \sums{d | q \\ (d, r) = 1} \mu(d) \sums{j \leq q/d} e_{rq}(a(bd^{\phi(r)} + j rd)) \\
&=& \sums{d | q \\ (d, r) = 1} \mu(d) e_{rq}(abd^{\phi(r)}) \sums{j \leq q/d} e_{q/d}(aj). 
\end{eqnarray*}
The last sum is non-zero only for $d=q$ and the claim follows.
\end{proof} 

We are now ready to prove the following approximation lemma for the generating function. 

\begin{Lemma} \label{lemma_generating_function}  Let $a, q \in \N$ with $1 \leq a \leq q$ and $(a, q) = 1$. Let $\alpha \in [0, 1]$, $\beta := \alpha - a/q$ and $N, b, r \in \N$ be such that $(b, r) = 1$.
Then
\begin{displaymath}
S_{b, r}(N, \alpha) = \frac{\mathbf{1}_{(q, r) = 1} \mu(q)e_{rq}(ab(q^{\phi(r)} - 1))r}{\phi(rq)} \sums{n \leq N} e(\beta n) + O\Big(q(1 + |\beta| N)E_{rq}(rN + b) \Big).
\end{displaymath}
\end{Lemma}

\begin{proof} 
Recalling the definition $(\ref{definition_of_E})$ and rearranging we see that
\begin{eqnarray*}
S_{b, r}(t, a/q)
& =& \sums{n \leq rt + b\\ n \equiv b \Mod r} \Lambda(n) e_{rq}(a(n - b))  \\
&=& e_{rq}(-a b)\sums{n \leq rt + b\\ n \equiv b \Mod r} \Lambda(n) e_{rq}(a n)\\
&=& e_{rq}(-a b) \sums{h = 1\\ h \equiv b \Mod r\\ (h, q) = 1}^{rq} e_{rq}(ah) \sums{n \leq rt + b\\ n \equiv h\Mod{rq}} \Lambda(n) + O(q)\\
&=& e_{rq}(-a b)\sums{h = 1\\ h \equiv b \Mod r\\ (h, q) = 1}^{rq} e_{rq}(ah) \Big(\frac{rt}{\phi(rq)} + O(E_{rq}(rt+b)) \Big) + O(q)\\
&=& e_{rq}(-a b) \frac{rt}{\phi(rq)}\sums{h = 1\\ h \equiv b \Mod r\\ (h, q) = 1}^{rq} e_{rq}(ah) + O(qE_{rq}(rt+b)) 
\end{eqnarray*}
Using Lemma \ref{lemma_rational_exp_sum} we now have
\begin{displaymath}
S_{b, r}(t, a/q) = \mathbf{1}_{(q, r) = 1} \mu(q)e_{rq}(ab(q^{\phi(r)} - 1))\frac{ rt}{\phi(rq)} + O(qE_{rq}(rt + b)).
\end{displaymath}
Thus 
\begin{displaymath}
U(t) := S_{b, r}(t, a/q) -  \frac{\mathbf{1}_{(q, r) = 1} \mu(q)e_{rq}(ab(q^{\phi(r)} - 1))r}{\phi(rq)}  \sums{n \leq t} 1 \ll qE_{rq}(rt+b).
\end{displaymath}
By partial summation it follows that 
\begin{eqnarray*}
&& S_{b, r}(N, \alpha) - \frac{\mathbf{1}_{(q, r) = 1} \mu(q)e_{rq}(ab(q^{\phi(r)} - 1))r}{\phi(rq)} \sums{n \leq N} e(\beta n) \\ 
&=& \sums{n \leq N} \Big(\Lambda(rn + b) e_q(a n)  - \frac{\mathbf{1}_{(q, r) = 1} \mu(q)e_{rq}(ab(q^{\phi(r)} - 1))r}{\phi(rq)}\Big) e(\beta n) \\
&=& U(N) e(\beta N) - \int_1^N U(t) 2\pi i \beta e(\beta t) dt \\
&\ll & qE_{rq}(rN+b) + q E_{rq}(rN+b)  N |\beta|. \qedhere
\end{eqnarray*}
\end{proof}

If $q \leq \log^A N$ and $\alpha \in \mathfrak{M}_{A, N}(q, a)$, then the error term in the previous lemma is $O(E_{rq}(rN + b) \log^{2A} N)$.

\subsection{Main term}
We are now ready to prove the following proposition. 

\begin{Proposition} \label{proposition_major_arcs} Let $A, N> 0$. Let $r, b_1, b_2, H \in \N$ be such that  $b_1, b_2 \leq r$, $(r, b_1) = (r, b_2) = 1$ and $H \leq N$. Then
\begin{align*}
\Big| \int_{\mathfrak{M}_{A, N}} S_{b_1, r}(N, \alpha) S_{b_2, r}(N, \alpha) e(-\alpha H) d\alpha  - & \frac{r^2}{\phi(r)^2}\mathfrak{S}_r(rH + b_1 +b_2)H \Big |  \\
& \ll \frac{\tau(Hr+b_1 +b_2)N}{\log^B N} + \sums{q \leq \log^A N}E_{rq}(rN + r) \log^{4A} N,
\end{align*}
for any $B > 0$, provided that $A$ is large enough depending on $B$.

\end{Proposition}

\begin{proof}

By Lemma \ref{lemma_generating_function} 
\begin{eqnarray*}
&& \int_{\mathfrak{M}_{A, N}} S_{b_1, r}(N, \alpha) S_{b_2, r}(N, \alpha) e(-\alpha H) d\alpha \\
&=& \sums{q \leq \log^A N\\ (q, r) = 1} \sums{1 \leq a \leq q \\ (a, q) = 1} \frac{\mu(q)e_{rq}(ab_1(q^{\phi(r)} - 1))r}{\phi(rq)} \frac{\mu(q)e_{rq}(ab_2(q^{\phi(r)} - 1))r}{\phi(rq)}e_q(-aH) \\
&& \times \int_{-\frac{\log^A N}{N}}^{\frac{\log^A N}{N}} \Big(\sums{n \leq N} e(\beta n) \Big)^2 e(-\beta H) d\beta\\
&& + O\Big( \sums{q \leq \log^A N} E_{rq}(rN + b_1) \log^{4A} N \Big).
\end{eqnarray*}
The sum in the last equation is called singular series and the integral is called singular integral.

By Lemma \ref{lemma_geometric_series}
\begin{eqnarray*}
\int_{ \frac{\log^A N}{N}}^{1/2} \Big(\sums{n \leq N} e(\beta n) \Big)^2 e(-\beta H) d\beta &\ll & \int_{ \frac{\log^A N}{N}}^{1/2} \beta^{-2} d\beta \\
& \ll & \frac{N}{\log^A N}. 
\end{eqnarray*}
Similarly 
\begin{displaymath}
\int_{-1/2}^{- \frac{\log^A N}{N}} \Big(\sums{n \leq N} e(\beta n) \Big)^2 e(-\beta H) d\beta  \ll \frac{N}{\log^A N}.
\end{displaymath}
Therefore the singular integral is 
\begin{eqnarray*}
 \int_{-\frac{\log^A N}{N}}^{\frac{\log^A N}{N}} \Big(\sums{n \leq N} e(\beta n) \Big)^2 e(-\beta H) d\beta 
 &=&  \int_{-1/2}^{1/2} \Big(\sums{n \leq N} e(\beta n) \Big)^2 e(-\beta H) d\beta + O\Big( \frac{N}{\log^A N}\Big) \\
 &=& H + O\Big( \frac{N}{\log^A N}\Big).
\end{eqnarray*}

Note that $e_q(-aH) = e_{rq}((q^{\phi(r)} - 1)raH)$. Hence the singular series equals
\begin{eqnarray*}
&& \sums{q \leq \log^A N\\ (q, r) = 1} \sums{1 \leq a \leq q \\ (a, q) = 1} \frac{\mu(q)e_{rq}(ab_1(q^{\phi(r)} - 1))r}{\phi(rq)} \frac{\mu(q)e_{rq}(ab_2(q^{\phi(r)} - 1))r}{\phi(rq)}e_q(-aH) \\
&=& \frac{r^2}{\phi(r)^2}\sums{q \leq \log^A N\\ (q, r) = 1} \frac{\mu(q)^2}{\phi(q)^2} \sums{1 \leq a \leq q\\ (a, q) = 1} e_q\Big(a\frac{q^{\phi(r)} - 1}{r}(b_1 + b_2 + rH) \Big) \\
&=& \frac{r^2}{\phi(r)^2}\sums{q \leq \log^A N\\ (q, r) = 1} \frac{\mu(q)^2}{\phi(q)^2}c_q(rH + b_1+b_2),
\end{eqnarray*}
where $c_q$ is the Ramanujan sum. We know (See e.g. \cite[(3.3)]{iwaniec-kowalski}) that
\begin{displaymath}
c_q(a) = \mu\Big( \frac{q}{(a, q)}\Big) \frac{\phi(q)}{\phi(q / (a, q))}.
\end{displaymath}
We see that
\begin{eqnarray*}
\sums{q > \log^A N\\ (q, r) = 1} \frac{\mu(q)^2}{\phi(q)^2} c_q(h) 
&=& \sums{k | h} \sums{q'k > \log^A N \\ (q', h/k) = 1} \frac{\mu(q'k)^2}{\phi(q'k)^2} c_{q'k}(h) \\
& \ll & \sums{k | h} \sums{q'k > \log^A N \\ (q', h/k) = 1} \frac{\mu(q'k)^2}{\phi(q'k)^2} \frac{\phi(q'k)}{\phi(q')} \\
& \ll & \sums{k | h} \frac{\mu(k)^2}{\phi(k)} \sums{q'k > \log^A N \\ (q', h/k) = 1} \frac{\mu(q')^2}{\phi(q')^2} \\
& \ll & \sums{k | h} \frac{\mu(k)^2}{\phi(k)} \sums{q' > \log^A N /k} \frac{1}{q'^{3/2}} \\
& \ll & \frac{\tau(h)}{\log^{B} N},
\end{eqnarray*}
for any $B > 0$, provided that $A$ is large enough depending on $B$. By Euler's product formula we have
\begin{displaymath}
\sums{q = 1\\ (q, r) = 1}^\infty \frac{\mu(q)^2}{\phi(q)^2} c_q(h) = \prod_{\substack{p \nmid r\\ p \nmid h}} \Big( 1 - \frac{1}{(p-1)^2}\Big) \prod_{\substack{p \nmid r\\ p | h}} \Big( 1 + \frac{1}{(p-1)}\Big) = \mathfrak{S}_r(h).
\end{displaymath}
Therefore
\begin{displaymath}
\sums{q \leq \log^A N\\ (q, r) = 1} \frac{\mu(q)^2}{\phi(q)^2}c_q(rH + b_1+b_2) =  \mathfrak{S}_r(rH + b_1+b_2) + O(\tau(rH + b_1 + b_2)\log^{-B} N).
\end{displaymath}
\end{proof}

\section{Minor arcs}\label{section_minor_arcs}
In this section, we aim to give satisfactory estimates for 
\begin{equation}\label{eq_minor_arcs_template_mean}
\max_{b_1: (b_1, r) = 1}\sums{n \leq N} \sums{b_2 \Mod r\\ (b_2, r) = 1}\Big| \int_{\mathfrak{m}_{A, N}} S_{b_1, r}(N, \beta)S_{b_2, r}(N, \beta)  e(-\beta n)  d\beta \Big|
\end{equation}
and for 
\begin{equation}\label{eq_minor_arcs_template}
\max_{b_1, b_2: (b_1b_2, r) = 1}\sums{n \leq N} \Big| \int_{\mathfrak{m}_{A, N}} S_{b_1, r}(N, \beta)S_{b_2, r}(N, \beta)  e(-\beta n)  d\beta \Big|.
\end{equation}

Let $\alpha \in [0, 1]$ and $a, q \in \N$ be such that $1 \leq a \leq q$, $(a, q) = 1$ and $|\alpha - a/q| \leq q^{-2}$. We assume this notation for the rest of the section. The key of giving sufficient estimates for (\ref{eq_minor_arcs_template_mean}) and (\ref{eq_minor_arcs_template}) is to study when 
\begin{equation} \label{eq_exp_sum_core}
\sums{n \leq X\\ n \equiv b \Mod r} \Lambda(n) e_r(\alpha n) = o\Big(\frac{X}{r} \Big).
\end{equation}
In order to analyse the previous sum, we will use Vaughan's identity for $\Lambda(n)$ and split the sum into type I and type II sums, which we will analyse in standard manners.

\subsection{Minor arcs with $b_2$-average}

In this subsection, we give an upper bound for (\ref{eq_minor_arcs_template_mean}). First, we introduce some auxiliary lemmas, starting with estimates for type I and type II sums.

\begin{Lemma} \label{lemma_type_I} \textbf{(Type I estimate)}  Let $M, N, X, A \geq 1$ with $MN \leq X$. Let $b, r \in \N$ with $(r, b) = 1$. Let $a_n$ be a complex sequence such that $|a_n| \leq 1$. Assume that $\log^{A+2} X \leq q \leq \frac{X}{r \log^{A+2} X}$ and $M \leq \frac{X}{r \log^{A+2} X}$ . Then

\begin{equation} \label{eq_type_I_sum}
\sums{mn \leq X\\ m \sim M, n \sim N \\ mn \equiv b \Mod r} a_m e_r(\alpha mn) \ll \frac{X}{r \log^A X}
\end{equation}
and
\begin{displaymath}
\sums{mn \leq X\\m \sim M, n \sim N\\ mn \equiv b \Mod r} a_m \log(n) e_r(\alpha mn) \ll \frac{X}{r \log^A X}.
\end{displaymath} 
\end{Lemma} 

\begin{proof}
Using Lemma \ref{lemma_geometric_series} and Lemma \ref{lemma_xy_exp_sum}  we see that

\begin{eqnarray*}
\sums{mn \leq X\\ m \sim M, n \sim N\\ mn \equiv b \Mod r} a_m e_r(\alpha mn)
 &\leq & \sums{m \sim M\\ (m, r) = 1} \Big| \sums{n \sim N\\ n \leq X/m\\ n \equiv b\overline{m} \Mod r} e_r(\alpha mn)\Big| \\
&\ll & \sums{m \sim M\\ (m, r) = 1} \max_{\substack{1 \leq a \leq r\\ (a, r) = 1}}\Big| \sums{rk + a \sim N\\ rk +  a \leq X/m} e(\alpha km)\Big| \\
&\ll & \sums{m \sim M} \min\Big(\frac{N}{r}, ||\alpha m||^{-1} \Big)\\
& \ll & (\frac{MN}{rq} + M + q) \log 2qM,
\end{eqnarray*}
which is $O(\frac{X}{r \log^A X})$, provided that 
\begin{displaymath}
 \log^{A+1} X \leq q \leq \frac{X}{r \log^{A+1} X} \text{ and } M \leq \frac{X}{r \log^{A+1} X}.
\end{displaymath}

Using partial summation we see that
\begin{eqnarray*}
\sums{mn \leq X\\m \sim M, n \sim N\\ mn \equiv b \Mod r} a_m \log(n) e_r(\alpha mn) 
&=& \log 2N \sums{mn \leq X\\m \sim M, n \sim N\\ mn \equiv b \Mod r} a_m e_r(\alpha mn) \\ 
&&  - \int_N^{2N} \sums{mn \leq X\\m \sim M, N \leq n < t\\ mn \equiv b \Mod r} a_m e_r(\alpha mn) \frac{dt}{t} \\
&\ll & \log N \max_{h \sim N} \Big| \sums{mn \leq X\\m \sim M, N \leq n < h\\ mn \equiv b \Mod r} a_m e_r(\alpha mn)\Big|.
\end{eqnarray*}
The remaining sum can be treated similarly as (\ref{eq_type_I_sum}).
\end{proof}

\begin{Lemma} \label{lemma_type_II} \textbf{(Type II estimate)} Let $M, N, X, A \geq 1$ with $MN \leq X$. Let $r \in \N$. Let $a_m, b_n$ be complex sequences such that $|a_m|, |b_n| \leq 1$. Assume that $\log^{A+1} X \leq q \leq \frac{X}{r \log^{A+1} X}$ and $\log^{A+1} X \leq M \leq \frac{X}{r \log^{A+1} X}$. Then
\begin{displaymath}
\sums{b \Mod r \\ (b, r) = 1} \Big|\sums{mn \leq X\\ m \sim M, n \sim N \\ mn \equiv b \Mod r} a_m b_n e_r(\alpha mn)\Big|^2  \ll \frac{X^2}{r \log^A X}.
\end{displaymath}

\end{Lemma}

\begin{proof}
Using Cauchy-Schwarz inequality and rearranging we see that

\allowdisplaybreaks
\begin{eqnarray*}
&& \sums{b \Mod r \\ (b, r) = 1} \Big|\sums{mn \leq X \\ m \sim M, n \sim N\\ mn \equiv b \Mod r} a_m b_n e_r(\alpha mn)\Big|^2 \\
&\leq & M \sums{b \Mod r \\ (b, r) = 1} \sums{m \sim M\\ (m, r) = 1}\Big|\sums{ n \sim N \\ n \leq X /m \\ n \equiv b\overline{m} \Mod r} b_n e_r(\alpha mn)\Big|^2 \\
&\leq & M \sums{m \sim M\\ (m, r) = 1} \sums{c \Mod r \\ (c, r) = 1} \Big|\sums{ n \sim N\\ n \leq X/m \\ n \equiv c \Mod r} b_n e_r(\alpha mn)\Big|^2 \\
&\leq & M   \sums{ m \sim M} \sums{c \Mod r \\ (c, r) = 1}\sums{n \sim N \\n \leq X /m \\ n \equiv c \Mod r}  \sums{n' \sim N\\ n' \leq X/m \\ n \equiv n' \Mod r} b_n \overline{b_{n'}} e_r(\alpha m (n - n')) \\
&\leq & M  \sums{n \sim N \\ (n, r) = 1} \sums{n' \sim N\\ n \equiv n' \Mod r} \Big| \sums{m \sim M\\ m \leq \min(X/n, X/n')} e_r(\alpha m (n - n')) \Big|.
\end{eqnarray*}

Using Lemma \ref{lemma_geometric_series} and Lemma \ref{lemma_xy_exp_sum} we have
\begin{eqnarray*}
&& M\sums{n \sim N \\ (n, r) = 1} \sums{n' \sim N\\ n \equiv n' \Mod r} \Big| \sums{m \sim M\\ m \leq \min(X/n, X/n')} e_r(\alpha m (n - n')) \Big| \\
&\ll & M  \sums{n \sim N\\ (n, r) = 1} \sums{n' \sim N\\ n \equiv n' \Mod r} \min \Big(M, ||\alpha (n - n') / r||^{-1}\Big) \\
&\ll & M  \sums{n \sim N\\ (n, r) = 1} \sums{|k| \leq 2N/r} \min \Big(M, ||\alpha k||^{-1}\Big) \\
&\ll & M  \sums{n \sim N\\ (n, r) = 1} \sums{1 \leq k \leq 2N/r} \min \Big(M, ||\alpha k||^{-1}\Big) + M^2N\\
&\ll &  MN \Big(\frac{MN}{rq} + \frac{N}{r} + q \Big) \log 2qN + M^2N,
\end{eqnarray*}
which is $O(\frac{X^2}{ r \log^A X})$, provided that
\begin{displaymath}
\log^{A+1} X \leq q \leq \frac{X}{r\log^{A+1} X} \text{ and } \log^{A+1} X \leq  M \leq \frac{X}{r\log^{A} X}.\qedhere
\end{displaymath} 
\end{proof}

Using Vaughan's identity and Lemmas \ref{lemma_type_I} and \ref{lemma_type_II}, we now can prove the following crucial lemma. 

\begin{Lemma} \label{lemma_mean_value_minor_arcs} Let $A, B \geq 1$, $r \in \N$ and $\log^{A} X \leq q \leq \frac{X}{r \log^{A} X}$. Assume that $r \leq X^{1/2} / \log^{A} X$ and $A$ is large enough depending on $B$. Then
\begin{displaymath}
\sums{b \Mod r \\ (b, r) = 1} \Big|\sums{n \leq X\\ n \equiv b \Mod r} \Lambda(n) e_r(\alpha n)\Big|^2 \ll \frac{X^2}{r \log^B X}.
\end{displaymath}
\end{Lemma}

\begin{proof} 
By Vaughan's identity and Cauchy-Schwarz inequality we get that, for any $y \geq 1$,
\allowdisplaybreaks
\begin{eqnarray*}
&& \sums{b \Mod r \\ (b, r) = 1} \Big|\sums{n \leq X \\ n \equiv b \Mod r} \Lambda(n) e_r(\alpha n)\Big|^2 \\
& \ll &  \sums{b \Mod r \\ (b, r) = 1} \Big|\sums{mn \leq X\\ mn \equiv b \Mod r\\ m \leq y} \mu(m) \log(n) e_r(\alpha mn)\Big|^2 \\
&& +  \sums{b \Mod r \\ (b, r) = 1} \Big|\sums{mn \leq X\\ mn \equiv b \Mod r\\ m \leq y^2} \Big(\sums{kl = m\\ k, l \leq y} \mu(k) \Lambda(l)\Big) e_r(\alpha mn)\Big|^2   \\
&& +  \sums{b \Mod r \\ (b, r) = 1} \Big|\sums{mn \leq X\\ mn \equiv b \Mod r\\ m > y} \mu(m) \Big(\sums{k | n\\ k > y} \Lambda(k)\Big) e_r(\alpha mn)\Big|^2   \\
& \ll &  (\log X)^2\sums{b \Mod r \\ (b, r) = 1} \sums{2^i 2^j \leq X\\ 2^i \leq y} \Big|\sums{mn \leq X\\ m \sim 2^i, n \sim 2^j\\ mn \equiv b \Mod r} \Big(\mu(m)\mathbf{1}_{m \leq y} \Big) \log(n) e_r(\alpha mn)\Big|^2 \\
&& +  (\log X)^2 \sums{b \Mod r \\ (b, r) = 1}  \sums{2^i 2^j \leq X \\2^i \leq y^2} \Big|\sums{mn \leq X\\m \sim 2^i, n \sim 2^j \\ mn \equiv b \Mod r} \mathbf{1}_{ m \leq y^2}\Big(\sums{kl = m\\ k, l \leq y} \mu(k) \Lambda(l)\Big) e_r(\alpha mn)\Big|^2 \\
&& + (\log X)^2 \sums{2^i 2^j \leq X\\ 2^i, 2^j > y/2} \sums{b \Mod r \\ (b, r) = 1}  \Big|\sums{mn \leq X\\ m \sim 2^i, n \sim 2^j \\ mn \equiv b \Mod r} \Big(\mu(m)\mathbf{1}_{m > y} \Big) \Big(\sums{k | n\\ k > y} \Lambda(k)\Big) e_r(\alpha mn)\Big|^2.
\end{eqnarray*}
Let $y = X^{1/5}$. Note that $\sums{h | n} \Lambda(h) = \log n$. Hence we can deal the first sum and the second sum with Lemma \ref{lemma_type_I}. Since  $X^{1/2} \leq \frac{X}{r \log^{A} X}$, we can deal the third sum with Lemma \ref{lemma_type_II}. 
\end{proof}

Using Lemma \ref{lemma_mean_value_minor_arcs} we can now give a sufficient upper bound for (\ref{eq_minor_arcs_template_mean}).

\begin{Proposition} \label{proposition_minor_arcs} Let $A \geq 1$. Let $N, r \in \N$ be such that $r \leq (Nr)^{1/2} / \log^A N$. Then

\begin{displaymath}
\max_{b_1: (b_1, r) = 1} \sums{n \leq N}  \sums{b_2 \Mod r\\ (b_2, r) = 1} \Big| \int_{\mathfrak{m}_{A, N}} S_{b_1, r}(N, \beta)S_{b_2, r}(N, \beta)  e(-\beta n)  d\alpha \Big| \ll \frac{N^2 r}{ \log^D N},
\end{displaymath}
for any $D > 0$, provided that $A$ is large enough depending on $D$.

\end{Proposition}

\begin{proof} Fix any $b_1 \Mod r$ with $(b_1, r) = 1$. Let 
\begin{displaymath}
f_{b_2}(\beta) := \mathbf{1}_{\beta \in \mathfrak{m}_{A, N}} S_{b_1, r}(N, \beta)S_{b_2, r}(N, \beta)
\end{displaymath}
and
\begin{displaymath}
\widehat{f}_{b_2}(n) := \int_0^1 f(\beta)e(-\beta n )d \beta.
\end{displaymath}
Using Cauchy-Schwarz inequality and Parseval's identity we get that
\begin{eqnarray*}
&& \Big(\sums{n \leq N} \sums{b_2 \Mod r\\ (b_2, r) = 1} \Big| \int_{\mathfrak{m}_{A, N}} S_{b_1, r}(N, \beta)S_{b_2, r}(N, \beta)  e(-\beta n)  d\beta \Big| \Big)^2 \\
&\leq & N r \sums{n \leq N} \sums{b_2 \Mod r\\ (b_2, r) = 1} \Big| \int_{\mathfrak{m}_{A, N}} S_{b_1, r}(N, \beta)S_{b_2, r}(N, \beta)  e(-\beta n)  d\beta \Big|^2 \\
&=& Nr \sums{b_2 \Mod r\\ (b_2, r) = 1} \sums{n \leq N} |\widehat{f}(n)|^2 \\
&=& Nr \sums{b_2 \Mod r\\ (b_2, r) = 1} \int_{\mathfrak{m}_{A, N}} |S_{b_1, r}(N, \beta)|^2 |S_{b_2, r}(N, \beta)|^2 d\beta.
\end{eqnarray*}
Let $\gamma \in \mathfrak{m}_{A, N}$. By Dirichlet's theorem (See e.g. \cite[Theorem 4.1]{nathanson}) there exist $a, q \in \N$ with $q \leq N / \log^A N =: T$, $1 \leq a \leq q$ and $(a, q) = 1$ such that $|\gamma - a/q| \leq (qT)^{-1} \leq q^{-2}$. Since $\gamma \in \mathfrak{m}_{A, N}$, we have that $q \geq \log^A N$. Hence by Lemma \ref{lemma_mean_value_minor_arcs} it follows that
\begin{displaymath}
\sums{b_2 \Mod r\\ (b_2, r) = 1}  | S_{b_2, r}(N, \gamma)|^2 \ll_\epsilon \frac{r N^2}{\log^C N}, 
\end{displaymath}
for any $C > 0$, provided that $A$ is large enough depending on $C$. We also see that
\begin{displaymath}
\int_{0}^1  |S_{b_1, r}(N, \beta)|^2 d\beta = \sums{n \leq N} \Lambda(rn + b_1)^2 \ll N \log^2 N.
\end{displaymath}
Therefore 
\begin{eqnarray*}
&& Nr \sums{b_2 \Mod r\\ (b_2, r) = 1} \int_{\mathfrak{m}_{A, N}} |S_{b_1, r}(N, \beta)|^2 |S_{b_2, r}(N, \beta)|^2 d\beta \\
& \leq & Nr \max_{\gamma \in \mathfrak{m}_{A, N}}\sums{b_2 \Mod r\\ (b_2, r) = 1} |S_{b_2, r}(N, \gamma)|^2 \int_{0}^1  |S_{b_1, r}(N, \beta)|^2 d\beta \\
& \ll & \frac{N^4 r^2}{ \log^{C - 2} N}.
\end{eqnarray*}
By choosing $C$ large enough the claim follows.
\end{proof}

\subsection{Minor arcs without $b_2$-average}

In this subsection, we give an upper bound for (\ref{eq_minor_arcs_template}). First, we introduce some auxiliary lemmas.

\begin{Lemma} \label{lemma_type_II_second} \textbf{(Alternative Type II estimate)} Let $M, N, X, A \geq 1$ with $MN \leq X$. Let $b, r \in \N$ with $(r, b) = 1$. Let $a_m, b_n$ be complex sequences such that $|a_m|, |b_n| \leq 1$. Assume that $\log^{A+1} X \leq q / (q, r) \leq \frac{X}{r^2 \log^{A+1} X}$ and $r \log^{A+1} X \leq M \leq \frac{X}{r \log^{A+1} X}$. Then
\begin{displaymath}
\sums{mn \leq X\\ m \sim M, n \sim N \\ mn \equiv b \Mod r} a_m b_n e_{rq}(a mn)  \ll \frac{X}{r \log^{A / 2} X}.
\end{displaymath}
\end{Lemma}

\begin{proof}
First we see that 
\begin{align*}
\sums{mn \leq X \\ m \sim M, n \sim N\\ mn \equiv b \Mod r} a_m b_n e_{rq}(a mn) 
&= \sums{i, j \Mod r \\ ij \equiv b \Mod r} \sums{mn \leq X \\ m \sim M, n \sim N\\ m \equiv i \Mod r  \\ n \equiv j \Mod r} a_m b_n e_{rq}(a mn) \\
& \ll r \max_{i, j} \Big| \sums{mn \leq X \\ m \sim M, n \sim N\\ m \equiv i \Mod r  \\ n \equiv j \Mod r} a_m b_n e_{rq}(amn) \Big|.
\end{align*}

We argue quite similarly as in the proof of Lemma \ref{lemma_type_II}. Using Cauchy-Schwarz inequality, Lemma \ref{lemma_geometric_series} and Lemma \ref{lemma_xy_exp_sum} we see that

\allowdisplaybreaks
\begin{align*}
& \Big|\sums{mn \leq X \\ m \sim M, n \sim N\\ m \equiv i \Mod r  \\ n \equiv j \Mod r} a_m b_n e_{rq}(a mn)\Big|^2 \\
&\leq  \frac{M}{r} \sums{m \sim M \\ m \equiv i \Mod r}\Big|\sums{ n \sim N \\ n \leq X /m \\ n \equiv j \Mod r} b_n e_{rq}(amn)\Big|^2 \\
& =  \frac{M}{r} \sums{m \sim M \\ m \equiv i \Mod r} \sums{n \sim N \\n \leq X /m \\ n \equiv j \Mod r}  \sums{n' \sim N\\ n' \leq X/m \\ n \equiv n' \Mod r} b_n \overline{b_{n'}} e_{rq}(a m (n - n')) \\
& \leq \frac{M}{r} \sums{n \sim N \\ n \equiv j \Mod r}  \sums{n' \sim N\\ n \equiv n' \Mod r} \Big| \sums{m \sim M\\ m \leq \min(X/n, X/n')\\ m \equiv i \Mod r} e_{rq}(a m (n - n')) \Big| \\
&\ll  \frac{M}{r} \sums{n \sim N \\ n \equiv j \Mod r}  \sums{n' \sim N\\ n \equiv n' \Mod r} \min \Big(M / r, ||a (n - n') / q||^{-1}\Big) \\
&\ll  \frac{MN}{r^2}  \sums{|k| \leq 2N/r} \min \Big(M / r, ||a k r / q||^{-1}\Big) \\
&\ll  \frac{MN}{r^2} \Big(\frac{MN}{r^2 q / (q, r)}+ \frac{N}{r} + \frac{q}{(q, r)}\Big)\log 2Nq  + \frac{M^2N}{r^3},
\end{align*}
which is $O(\frac{X^2}{ r^4 \log^A X})$, provided that
\begin{displaymath}
\log^{A+1} X \leq q / (q, r) \leq \frac{X}{r^2\log^{A+1} X} \text{ and } r \log^{A+1} X \leq  M \leq \frac{X}{r\log^{A} X}.
\end{displaymath}
\end{proof}

\begin{Lemma} \label{lemma_rational_estimate} Let $A, B, C \geq 0$ and $b, r \in \N$ with $(r, b) = 1$. Let $\beta \in \R$ be such that $|\beta| \leq X^{-1}r \log^B X$. Assume that $r \leq X^{1/3} / \log^{A} X$, $\log^{A} X \leq q / (q, r) \leq X/ (r^2 \log^{A} X)$ and $A$ is large enough depending on $B$ and $C$. Then
\begin{displaymath}
\sums{n \leq X \\ n \equiv b \Mod r} \Lambda(n) e_{r}\Big(\Big(\frac{a}{q} + \beta \Big)n\Big) \ll \frac{X}{r \log^C X}.
\end{displaymath}
\end{Lemma}

\begin{proof} Using partial summation, we see that 
\begin{align*}
\sums{n \leq X \\ n \equiv b \Mod r} \Lambda(n) e_{rq}(an)e_r(\beta n)
= & \ e_r(\beta X) \sums{n \leq X \\ n \equiv b \Mod r} \Lambda(n) e_{rq}(an) \\
& - \int_1^X \frac{2 \pi i \beta}{r}\sums{n \leq t \\ n \equiv b \Mod r} \Lambda(n) e_{rq}(an) e_r(\beta t) dt \\
 \ll & \ \log^B X \max_{t \leq X} \Big|\sums{n \leq t \\ n \equiv b \Mod r} \Lambda(n) e_{rq}(an) \Big|.
\end{align*}
By Vaughan's identity 
\allowdisplaybreaks
\begin{align*}
& \sums{n \leq t \\ n \equiv b \Mod r} \Lambda(n) e_{rq}(an)\\
= & \ \sums{mn \leq t\\ mn \equiv b \Mod r\\ m \leq y} \mu(m) \log(n) e_{rq}(a mn)
+  \sums{mn \leq t\\ mn \equiv b \Mod r\\ m \leq y^2} \Big(\sums{kl = m\\ k, l \leq y} \mu(k) \Lambda(l)\Big) e_{rq}(a mn)   \\
& +  \sums{mn \leq t\\ mn \equiv b \Mod r\\ m > y} \mu(m) \Big(\sums{k | n\\ k > y} \Lambda(k)\Big) e_{rq}(a mn)   \\
 = & \  \sums{2^i 2^j \leq t\\ 2^i \leq y} \sums{mn \leq t\\ m \sim 2^i, n \sim 2^j\\ mn \equiv b \Mod r} \Big(\mu(m)\mathbf{1}_{m \leq y} \Big) \log(n) e_{rq}(a mn) \\
& +  \sums{2^i 2^j \leq t \\2^i \leq y^2} \sums{mn \leq t\\m \sim 2^i, n \sim 2^j \\ mn \equiv b \Mod r} \mathbf{1}_{ m \leq y^2}\Big(\sums{kl = m\\ k, l \leq y} \mu(k) \Lambda(l)\Big) e_{rq}(a mn) \\
& + \sums{2^i 2^j \leq t\\ 2^i, 2^j > y / 2} \sums{mn \leq t\\ m \sim 2^i, n \sim 2^j \\ mn \equiv b \Mod r} \Big(\mu(m)\mathbf{1}_{m > y} \Big) \Big(\sums{k | n\\ k > y} \Lambda(k)\Big) e_{rq}(a mn).
\end{align*}
Let $y = r \log^A X$. By assumption we have $y^2 \leq X/(r \log^A X)$. Note that $\sums{h | n} \Lambda(h) = \log n$. Hence we can deal the first sum and the second sum with Lemma \ref{lemma_type_I}. We deal the third sum with Lemma \ref{lemma_type_II_second}. 
\end{proof}

Since we have requirement $q /(q, r) \leq X / (r^2 \log^{A + 1} X)$ in Lemma \ref{lemma_type_II_second}, estimate $(\ref{eq_exp_sum_core})$ does not hold for entire minor arcs. Therefore we need the following lemma, which follows from the proof of \cite[Theorem 2]{bauer_old} (Particularly, the condition $|\alpha - a'/q'| \leq 1/q'$ comes from \cite[Lemma 2.1]{bauer_old}).

\begin{Lemma} \label{lemma_bauer_th2_original} Let $X \geq 1$, $\alpha' \in [0, 1]$ and $a', q' \in \N$ be such that $|\alpha - a'/q'| \leq 1/q'$. Let $b, r \in \N$ be such that $(b, r) = 1$. Write $\beta' := \alpha' - a'/q'$.  Then
\begin{align*}
& \sums{n \leq X \\ n \equiv b \Mod r} \Lambda(n) e(\alpha' n) \\
& \ll (\log X )^{O(1)} \tau([q', r])^{1/2} \frac{(q' / h)^{1/2}}{[q', r]} \Big( [q', r] X^{1/2}\sqrt{1 + |\beta'| X} + ([q', r])^{1/2} X^{4/5} + \frac{X}{\sqrt{1 + |\beta'| X } }\Big) + q' \log X,
\end{align*}
where 
\begin{displaymath}
h = \prod_{\substack{p\\ p^t || (q, r)\\ p^t || q}} p^t.
\end{displaymath}
\end{Lemma}

We rewrite the previous lemma in the following more suitable form.

\begin{Lemma} \label{lemma_bauer_th2} Let $X \geq 1$. Let $b, r \in \N$ be such that $(b, r) = (q, r) = 1$. Write $\beta := \alpha - a/q$ and assume that $|\beta| \leq 1/q$.  Then
\begin{align*}
& \sums{n \leq X \\ n \equiv b \Mod r} \Lambda(n) e_r(\alpha n) \\
& \ll (\log X )^{O(1)} \tau(qr) \frac{q^{1/2}}{qr} \Big( qr X^{1/2}\sqrt{1 + |\beta| X / r} + (qr)^{1/2} X^{4/5} + \frac{X}{\sqrt{1 + |\beta| X / r} }\Big) + qr \log X.
\end{align*}
\end{Lemma}

\begin{proof} Set $q' = q r / (a, r)$, $a' = a / (a, r)$ and $\alpha' = \alpha /r$. The claim now follows from Lemma \ref{lemma_bauer_th2_original} as $q' / h = q$,  $[q', r] = qr$ and $\beta' = \beta / r$.
\end{proof}

Using the previous lemmas we now prove the following important lemma. 

\begin{Lemma} \label{lemma_minor_arcs_no_mean_estimate} Let $\epsilon > 0$, $A, C > 0$, $X > 1$ and $\alpha \in \mathfrak{m}_{A+1, X/r}$. Let $b, r \in \N$ with $(b, r) = 1$. Assume that $r \leq X^{1/3 - \epsilon}$. Then
\begin{displaymath}
\sums{n \leq X \\ n \equiv b \Mod r} \Lambda(n) e_r(\alpha n) \ll_\epsilon \frac{\tau(r)X}{r \log^C X} + \max_{\substack{q < r \log^A X \\ (r, q) \neq 1}} r E_{rq}(X) \log^A X,
\end{displaymath}
provided that $A$ is large enough depending on $C$.
\end{Lemma}

\begin{proof} By Dirichlet's theorem there exist $q \leq X / (r^2 \log^A X)$, $1 \leq a \leq q$ such that $(a, q) = 1$ and $|\alpha - a/q| \leq \frac{r^2 \log^A X}{qX}$. Write $\beta := \alpha - a/q$. Let $B > 0$ to be chosen later. We split into the following cases
\begin{enumerate}
\item[\textbf{I:}] $q \leq \log^A X$ and $|\beta| \leq \frac{r\log^B X}{X}$.
\item[\textbf{II:}] $q > \log^A X$, $(q, r) = 1$ and $|\beta| \leq \frac{r\log^B X}{X}$.
\item[\textbf{III:}] $q \geq r \log^A X $. 
\item[\textbf{IV:}] $q < r \log^A X$, $(q, r) = 1$ and $|\beta| > \frac{r\log^B X}{X}$. 
\item[\textbf{V:}] $q < r \log^A X$ and $(r, q) \neq 1$.
\end{enumerate}

\textbf{Case I: } By definition we have $\alpha \in \mathfrak{M}_{A+1, X/r}$, provided that $B < A$. So this case cannot actually occur.

\textbf{Case II: } Follows from Lemma \ref{lemma_rational_estimate}, provided that $A$ is large enough depending on $B$ and $C$. 

\textbf{Case III: } Follows from Lemma \ref{lemma_rational_estimate}, provided that $A$ is large enough depending on $C$.

\textbf{Case IV: } Using Lemma \ref{lemma_bauer_th2} we see that 
\begin{align*}
\sums{n \leq X \\ n \equiv b \Mod r} \Lambda(n) e_r(\alpha n) \ll_\epsilon (\log^{O(1)} X) \tau(r) \Big( X^{1/2 + \epsilon / 2}r^{1/2} + \frac{X^{4/5 + \epsilon}}{r^{1/2}}  + \frac{X}{r \log^{B/2} X}\Big) + r^2 \log^{A + 1} X,
\end{align*}
which is $O(\tau(r) X /(r \log^C X))$, provided that $B$ is large enough depending on $C$ and $r \leq X^{1/3 - \epsilon}$.

\textbf{Case V: } By Lemma \ref{lemma_generating_function}
\begin{displaymath}
\sums{n \leq X \\ n \equiv b \Mod r} \Lambda(n) e_r(\alpha n) \ll r E_{rq}(X) \log^A X.
\end{displaymath} \qedhere
\end{proof}

Using Lemma \ref{lemma_minor_arcs_no_mean_estimate} we can now give sufficient upper bound for (\ref{eq_minor_arcs_template}).

\begin{Proposition} \label{proposition_minor_arcs_second}
Let $A \geq 1$ and $\epsilon > 0$. Let $N, r \in \N$ be such that $r \leq (Nr)^{1/3 - \epsilon}$. Then
\begin{align*}
\max_{b_1, b_2: (b_1b_2, r) = 1}  & \sums{n \leq N} \Big| \int_{\mathfrak{m}_{A, N}} S_{b_1, r}(N, \beta)S_{b_2, r}(N, \beta)  e(-\beta n)  d\beta \Big| \\ 
& \ll_\epsilon \frac{N^2 \tau(r)}{ \log^B N} + N^{3/2} \Big(r\max_{\substack{q < r\log^A N \\ (r, q) \neq 1}} E_{r q}(rN + r) \Big)^{1/2} \log^{A/2 + 2} N ,
\end{align*}
for any $B > 0$, provided that $A$ is large enough depending on $B$.
\end{Proposition}

\begin{proof} Arguing as in Proposition \ref{proposition_minor_arcs}.
\begin{align*}
\Big(\max_{b_1, b_2: (b_1b_2, r) = 1}  \sums{n \leq N} \Big| \int_{\mathfrak{m}_{A, N}} S_{b_1, r}(N, \beta)S_{b_2, r}(N, \beta)  e(-\beta n)  d\beta \Big|\Big)^2
& \ll N^2 \log^2 N \max_{\gamma \in \mathfrak{m}_{A, N}} | S_{b_1, r}(N, \gamma) |^2 \\
& \ll N^3 \log^2 N \max_{\gamma \in \mathfrak{m}_{A, N}} | S_{b_1, r}(N, \gamma) |
\end{align*}
and the claim follows from Lemma  \ref{lemma_minor_arcs_no_mean_estimate}.
\end{proof}

\section{Error term} \label{section_error_term}

In this section, we prove two error term estimates. The first one follows from Bombieri-Vinogradov theorem.

\begin{Lemma}\label{lemma_bombieri_vinogradov} Let $A, B, N > 0$. Let $q \in \N$ with $q \leq \log^A N$. Let $R \leq (RN)^{1/2} / \log^{A + B + 5} N$. Then
\begin{displaymath}
\sums{r \leq R}  |E_{qr}(rN)| \ll \frac{RN}{\log^B N}.
\end{displaymath}
\end{Lemma}

\begin{proof}
By the Bombieri-Vinogradov theorem (See e.g. \cite[Theorem 8.34]{tenenbaum})
\begin{displaymath}
\sums{r \leq R}  |E_{qr}(rN )| \leq \sums{r \leq R}  |E_{qr}(RN )| \leq \sums{r \leq R \log^A N}  |E_{r}(RN )|  \leq \sums{r \leq (RN)^{1/2} / \log^{B+5} N}  |E_{r}(RN)| \ll \frac{RN}{\log^{B} N}.
\end{displaymath}
\end{proof}

For the proof of the second error term estimate we need information about zeros of Dirichlet $L$-functions $L(s, \chi)$. By \cite[Theorem 5.24]{iwaniec-kowalski} and \cite[Theorem 10.4]{iwaniec-kowalski} we have the following density lemma. 
\begin{Lemma} \label{lemma_density_theorem} Let $1/2 \leq \alpha \leq 1$ and $\epsilon > 0$. Let $r \in \N$. There exists $D = D(\epsilon)> 0$ such that
\begin{displaymath}
\sums{\psi \Mod r} \sums{L(\rho, \psi) = 0\\ |Im(\rho)| \leq H\\ Re(\rho) \geq \alpha} 1 \ll_\epsilon \Big((rH)^{(2 + \epsilon)(1 - \alpha)} + (rH)^{c(\alpha)(1 - \alpha)} \Big) \log^D (rH),
\end{displaymath}
where
\begin{displaymath}
c(\alpha) = \min \Big( \frac{3}{2 - \alpha}, \frac{3}{3\alpha - 1} \Big).
\end{displaymath}

\end{Lemma}

Let $E > 0$ and write $K := \frac{E\log \log N}{\log N}$. Define 
\begin{align*}
\mathcal{D}_E &:= \{s \in \C :  Re(s) > 1 - K \text{ and } |Im(s)| \leq N\}, \\
Z_{E}(N) & := \{h \leq N \mid \exists s \in \mathcal{D}_E, \text{ primitive }  \chi \Mod{h}: L(s, \chi) = 0 \}, \\
Z_E(N, Q) &:= \{h \leq N/Q \mid \exists q \leq Q: hq \in Z_E(N)\}.
\end{align*}

\begin{Lemma} \label{lemma_size_of_Z} Let $E > 0$ and $Q, N > 1$. Then,
\begin{equation*} 
|Z_E(N, Q)| \ll_\epsilon \min\Big(N^\epsilon, Q\log^D N\Big),
\end{equation*}
for any $\epsilon > 0$ and for some $D > 0$ depending on $E$.
\end{Lemma}

\begin{proof}
First, using \cite[Theorem 10.4]{iwaniec-kowalski},  we see that 
\begin{align*}
|Z_E(N, Q)| 
& \leq  \sums{h \leq N/Q} \sums{q \leq Q} \sums{\chi \Mod{qh} \\ \chi \text{ primitive}} \sums{L(\rho, \chi) = 0\\ \rho \in \mathcal{D}_E} 1 \\
& \leq \sums{q \leq Q} \sums{h \leq N} \sums{\chi \Mod{h} \\ \chi \text{ primitive}} \sums{L(\rho, \chi) = 0\\ \rho \in \mathcal{D}_E} 1 \\
& \ll Q \log^D N,
\end{align*}
for some $D > 0$ depending on $E$. Similarly, we have
\begin{align*}
|Z_E(N, Q)| 
& \leq  \sums{h \leq N/Q} \sums{q \leq Q} \sums{\chi \Mod{qh} \\ \chi \text{ primitive}} \sums{L(\rho, \chi) = 0\\ \rho \in \mathcal{D}_E} 1 \\
& \leq \sums{t \leq N} \tau(t) \sums{\chi \Mod{t} \\ \chi \text{ primitive}} \sums{L(\rho, \chi) = 0\\ \rho \in \mathcal{D}_E} 1 \\
& \leq  \max_{t \leq N}\tau(t) \sums{t \leq N}\sums{\chi \Mod{t} \\ \chi \text{ primitive}} \sums{L(\rho, \chi) = 0\\ \rho \in \mathcal{D}_E} 1 \\
& \ll_\epsilon N^\epsilon,
\end{align*}
for any $\epsilon > 0$.
\end{proof}

Using Lemma \ref{lemma_density_theorem}, we can now prove the second error term estimate. 

\begin{Lemma} \label{lemma_512} Let $E, \epsilon > 0$. Let $X \geq 3$ and $d \leq X^{5/12 - \epsilon}$ be such that, for all $h | d$, we have $h \not\in Z_E(X)$ \footnote{For $\zeta$-function $\mathcal{D}_E$ is zero-free region. Hence $1 \not\in Z_E(X)$.}. Then
\begin{displaymath}
E_{d}(X) \ll_\epsilon \frac{X \tau(d)}{\phi(d)\log^B X},
\end{displaymath}
for any $B > 0$, provided that $E$ is large enough depending on $B$ and $\epsilon$. 
\end{Lemma}

\begin{proof}
First we see that, for $(b, d) = 1$, 
\begin{eqnarray*}
\sums{n \leq X\\ n \equiv b \Mod{d}} \Lambda(n) 
&=& \frac{1}{\phi(d)}\sums{\chi \Mod{d}} \overline{\chi}(b)\sums{n \leq X} \Lambda(n)\chi(n) \\
&=& \frac{1}{\phi(d)} \sums{n \leq X\\ (n, d) = 1} \Lambda(n) + \frac{1}{\phi(d)} \sums{h | d \\h \neq 1} \sums{\chi \Mod{h} \\ \chi \text{ primitive}}  \overline{\chi}(b)\sums{n \leq X\\ (n, d) = 1}\Lambda(n)\chi(n).  \\
\end{eqnarray*}
By the prime number theorem $\sums{n \leq X\\ (n, d) = 1} \Lambda(n) = X + O(X \log^{-B} X)$, for any $B > 0$.  It remains to prove that 
\begin{displaymath}
\sums{\chi \Mod{h} \\ \chi \text{ primitive}} \Big| \sums{n \leq X} \Lambda(n)\chi(n) \Big| \ll \frac{X}{\log^B X},
\end{displaymath}
for any $B > 0$ and $h | d$ with $h \neq 1$. 

Using the explicit formula (see e.g. \cite[Proposition 5.25]{iwaniec-kowalski}) we see that
\begin{eqnarray*}
\sums{n \leq X} \Lambda(n)\chi(n) 
&=& \sums{L(\rho, \chi) = 0\\ |Im(\rho)| \leq T} \frac{X^\rho - 1}{\rho} + O \Big(\frac{X}{T} \log^2 X \Big)  \\
&\ll & \sums{2^i \leq T} \sums{L(\rho, \chi) = 0\\ |Im(\rho)| \sim 2^i \\ Re(\rho) \geq 1/2} \frac{X^{Re(\rho)}}{|Im(\rho)|} +  \sums{L(\rho, \chi) = 0\\ 0 \leq |Im(\rho)| \leq 1\\ Re(\rho) \geq 1/2} X^{Re(\rho)} + \frac{X}{T} \log^2 X  \\
&\ll & \sums{2^i \leq T} \frac{1}{2^i} \sums{L(\rho, \chi) = 0\\ |Im(\rho)| \leq 2^i\\ Re(\rho) \geq 1/2} X^{Re(\rho)} + \frac{X}{T} \log^2 X ,
\end{eqnarray*}
for any $T > 0$. We take $T = d^{1 + \eta}$ for $\eta > 0$ to be chosen later. Using Lemma \ref{lemma_density_theorem}, we see that, for any $U \geq 1$,

\allowdisplaybreaks
\begin{align*}
& \frac{1}{U}\sums{\chi \Mod{h} \\ \chi \text{ primitive}}  \sums{L(\rho, \chi) = 0\\ |Im(\rho)| \leq U  \\ Re(\rho) \geq 1/2} X^{Re(\rho)}  \\ 
& \ll \frac{1}{U}\sums{\chi \Mod{h} }   \sums{L(\rho, \chi) = 0\\ 1/2 \leq Re(\rho)  < 1 - K \\ |Im(\rho)| \leq U} X^{Re(\rho)} \\ 
& \ll \frac{1}{U} \log X \max_{\frac12 \leq \beta \leq 1 - K} X^{\beta + 1/ \log X} \sums{\chi \Mod{h} } \sums{L(\rho, \chi) = 0 \\ |Im(\rho)| \leq U \\ \beta \leq Re(\rho) \leq \beta +1/\log X}  1 \\
& \ll_\epsilon  \log^{D + 1} X \max_{1/2 \leq \beta \leq 1 - K} X^{\beta}\frac{1}{U}\Big( (dU)^{(2 + \epsilon')(1 - \beta)} + (dU)^{c(\beta)(1 - \beta)} \Big),
\end{align*}
for any $\epsilon' > 0$ and $D$ large enough depending on $\epsilon'$.  
Let $\lambda(\beta) = \max(2 + \epsilon', c(\beta))$. We can see that 
\begin{displaymath}
\lambda(\beta) = 
\begin{cases}
2 + \epsilon' & \text{ when } 1/2 \leq \beta \leq \frac{1 + 2 \epsilon'}{2 + \epsilon'}. \\
\frac{3}{2 - \beta} & \text{ when } \frac{1 + 2 \epsilon'}{2 + \epsilon'} \leq \beta \leq 3/4 \\
\frac{3}{3\beta - 1} & \text{ when } 3/4 \leq \beta \leq \frac{5 + \epsilon'}{6 + 3\epsilon'} \\
2 + \epsilon' & \text{ when } \frac{5 + \epsilon'}{6 + 3\epsilon'} \leq \beta \leq 1. \\
\end{cases}
\end{displaymath}
Write $d  := X^\theta$ and 
\begin{displaymath}
G(\beta) := X^{\beta}\frac{1}{U}(dU)^{\lambda(\beta)(1 - \beta)} .
\end{displaymath}
It then suffices to show that 
\begin{displaymath}
G(\beta) \ll \frac{X}{\log^B X},
\end{displaymath}
for any $B \geq 1$, whenever $1 \leq U \leq T$. We have the following cases:

\begin{enumerate}
\item[\textbf{I:}] $\lambda(\beta)(1 - \beta) \geq 1$.
\item[\textbf{II:}] $\lambda(\beta)(1 - \beta) < 1$.
\end{enumerate}

\textbf{Case I: } In this case $G(\beta)$ is maximal when $U = T = d^{1 + \eta}$. We also have that $\beta \leq \frac{1 + \epsilon'}{2 + \epsilon'}$ and so $\lambda(\beta) = 2 + \epsilon'$. Hence
\begin{displaymath}
G(\beta) \leq X^{\frac{1 + \epsilon'}{2 + \epsilon'} - \theta + \theta(2 + \eta)(2+ \epsilon')\frac12 } \ll X^{1 - \eta},
\end{displaymath}
provided that $\theta < 1/2$ and $\eta$ and $\epsilon'$ are sufficiently small depending on $\theta$. \\

\textbf{Case II: } In this case $G(\beta)$ is maximal when $U = 1$. We also see that $\lambda(\beta) \leq 12/5$. Since $\theta \leq \frac{5}{12} - \epsilon$ and $\beta \leq 1 - K$, we have
\begin{displaymath}
G(\beta) \leq X^{\beta + \theta \frac{12}{5}(1 - \beta)}  \leq X^{1 - \frac{12}{5}\epsilon K},
\end{displaymath} 
so the claim follows once $E$ is large enough depending on $B$ and $\epsilon$. 
\end{proof}

The following two lemmas are direct consequences of Lemma \ref{lemma_512}. 

\begin{Lemma} \label{lemma_512_1} Let $\epsilon > 0$, $A, E \geq 1$ and $N \geq 3$. Let $r$ be a prime such that $r \not\in Z_E(N, \log^A N)$ and $r \leq (rN)^{5/12 - \epsilon}$. Then
\begin{displaymath}
\max_{\substack{q < \log^A N}} E_{r q}(rN) \ll_\epsilon \frac{N}{ \log^B N},
\end{displaymath}
for any $B > 0$, provided that $E$ is large enough depending on $B$ and $\epsilon$.
\end{Lemma}

\begin{Lemma}\label{lemma_512_2} Let $\epsilon > 0$, $A, E \geq 1$ and $N \geq 3$. Let $r$ be a prime such that $r, r^2 \not\in Z_E(N, \log^A N)$ and $r^2 \leq (rN)^{5/12 - \epsilon}$. Then
\begin{displaymath}
\max_{\substack{q < r\log^A N \\ r | q}} E_{r q}(rN) \ll_\epsilon \frac{N}{r \log^B N},
\end{displaymath}
for any $B > 0$, provided that $E$ is large enough depending on $B$ and $\epsilon$.
\end{Lemma}

\noindent \textbf{Remark.} Assuming GRH it is known that $E_d(x) \ll \sqrt{x}\log x$ for any $d \leq x$.

\section{Proofs of theorems}\label{section_proofs}

In this section, we finish the proofs of Theorems \ref{theorem_technical_main}, \ref{theorem_technical_3} and \ref{theorem_technical_2}. First, we see that, for all $n \leq N$, 
\begin{align*}
&  \Big| \sums{n_1, n_2\\ n = n_1 + n_2} \Lambda(rn_1 + b_1)\Lambda(rn_2 + b_2) - \frac{r^2}{\phi(r)^2}\mathfrak{S}_r(rn + b_1 +b_2)n \Big | \\
& = \  \Big| \int_0^1 S_{b_1, r}(N, \alpha)S_{b_2, r}(N, \alpha)  e(-\alpha n)  d\alpha - \frac{r^2}{\phi(r)^2}\mathfrak{S}_r(rn + b_1 +b_2)n \Big | \\
& \leq \ \Big| \int_{\mathfrak{M}_{A, N}} S_{b_1, r}(N, \alpha)S_{b_2, r}(N, \alpha)  e(-\alpha n)  d\alpha - \frac{r^2}{\phi(r)^2}\mathfrak{S}_r(rn + b_1 +b_2)n \Big | \\
& \ \ \ + \  \Big| \int_{\mathfrak{m}_{A, N}} S_{b_1, r}(N, \alpha)S_{b_2, r}(N, \alpha)  e(-\alpha n)  d\alpha \Big |.
\end{align*}

By Proposition \ref{proposition_major_arcs}
\begin{align}
\Big| \int_{\mathfrak{M}_{A, N}} S_{b_1, r}(N, \alpha) S_{b_2, r}(N, \alpha) & e(-\alpha n) d\alpha  -   \frac{r^2}{\phi(r)^2}\mathfrak{S}_r(rn + b_1 +b_2)n \Big |  \nonumber\\
& \ll \frac{\tau(nr+b_1 +b_2)N}{\log^C N} + \sums{q \leq \log^A N}E_{rq}(rN + r) \log^{4A} N, \label{eq_major_arcs}
\end{align}
for any $C > 0$, provided that $A$ is large enough depending on $C$. By \cite{vinogradov-linnik}
\begin{equation} \label{eq_divisor_bound}
\sums{n \leq N} \tau(rn + b) \ll \frac{\phi(r)}{r} N \log N,
\end{equation}
for $r < N$ and $0 < b < r$ with $(b, r) = 1$. 

\begin{proof}[Proof of Theorem \ref{theorem_technical_main}] 
By (\ref{eq_major_arcs}), (\ref{eq_divisor_bound}) and Lemma \ref{lemma_bombieri_vinogradov}
\begin{align*}
& \sums{r \leq R} \max_{b_1: (b_1, r) = 1} \sums{b_2 \Mod r\\ (b_2, r) = 1} \sums{n \leq N} \Big| \int_{\mathfrak{M}_{A, N}} S_{b_1, r}(N, \alpha)S_{b_2, r}(N, \alpha)  e(-\alpha n)  d\alpha  - \frac{r^2}{\phi(r)^2}\mathfrak{S}_r(rn + b_1 +b_2)n \Big | \\
& \ll \ \frac{NR}{\log^C N}\sums{r \leq R} \max_{b_1, b_2 \Mod r} \sums{n \leq N} \tau(nr + b_1 + b_2) + NR\log^{4A} N \sums{q \leq \log^A N} \sums{r \leq R} E_{rq}(rN + r) \\
& \ll \ \frac{N^2 R^2}{\log^B N},
\end{align*}
for any $B$, provided that $A$ is large enough depending on $B$. The claim now follows from Proposition \ref{proposition_minor_arcs}.
\end{proof}

\begin{proof}[Proof of Theorem \ref{theorem_technical_3}] 
By (\ref{eq_major_arcs}), (\ref{eq_divisor_bound}) and Lemmas \ref{lemma_512_1} and \ref{lemma_size_of_Z} 
\begin{align*}
& \max_{b_1, b_2: (b_1b_2, r) = 1} \sums{n \leq N}  \Big| \int_{\mathfrak{M}_{A, N}} S_{b_1, r}(N, \alpha)S_{b_2, r}(N, \alpha)  e(-\alpha n)  d\alpha  - \frac{r^2}{\phi(r)^2}\mathfrak{S}_r(rn + b_1 +b_2)n \Big | \\
& \ll \ \frac{N}{\log^C N} \max_{b_1, b_2: (b_1b_2, r) = 1} \sums{n \leq N} \tau(nr + b_1 + b_2) + N(\log^{5A} N) \max_{q \leq \log^A N} E_{rq}(rN + r) \\
& \ll \frac{N^2}{\log^B N},
\end{align*}
for all but $O(\log^D N)$ primes $r$ and for any $B$, provided that $A$ is large enough depending on $B$ and $D$ is large enough depending on $A$ and $B$. 

By Proposition \ref{proposition_minor_arcs_second} we have that
\begin{align*}
\max_{b_1, b_2: (b_1b_2, r) = 1}  & \sums{n \leq N} \Big| \int_{\mathfrak{m}_{A, N}} S_{b_1, r}(N, \beta)S_{b_2, r}(N, \beta)  e(-\beta n)  d\beta \Big| \\ 
& \ll_\epsilon \frac{N^2 \tau(r)}{ \log^B N} + N^{3/2} \Big(r\max_{\substack{q < r\log^A N \\ (r, q) \neq 1}} E_{r q}(rN + r) \Big)^{1/2} \log^{A/2 + 2} N,
\end{align*}
for any $B > 0$, provided that $A$ is large enough depending on $B$. The claim follows using Lemmas \ref{lemma_512_2} and \ref{lemma_size_of_Z} to the second error term. 
\end{proof}

\begin{proof}[Proof of Theorem \ref{theorem_technical_2}] The proof is similar to the proof of Theorem \ref{theorem_technical_3}, but we use Proposition \ref{proposition_minor_arcs} instead of Proposition \ref{proposition_minor_arcs_second}.
\end{proof}

\subsection*{Acknowledgements}
The author wants to thank his supervisor Kaisa Matomäki for suggesting the topic and for many useful comments. During the work, the author was supported by Emil Aaltonen foundation.

\printbibliography

\end{document}